\documentclass[a4paper,11pt]{article}
\usepackage{setspace}
\usepackage{cite}
\usepackage{mathrsfs}
\usepackage{amsthm,amsmath,amssymb}

\title{Model Structures and Quantum Cohomology of Higher Orbifolds}
\author{Jiajun Dai}
\begin{document}
\maketitle
\section{Introduction}
This is an informal essay, though propped with detailed explanations in some parts, recording some thoughts of mine mainly about model structure on higher orbifolds and its application in quantum cohomology of higher orbifolds in recent months, developing Chen-Ruan's orbifold quantum cohomology. 
The thoughts originated from two ideas:\\

One is owing to the Taylor series expansion of the $ exponential \: map $ in Riemmanian geometry. 
Now that this formula which is so important in calculus can be transplanted to geometry, can it be generalized to other contexts? 
What will happen if sum became direct sum, and power became multiple (direct, tensor or other) product? $ Graded \: algebra $ like universal enveloping algebra, tensor algebra and exterior algebra provides a suitable answer, and graded category fits in some sense. 
Extending it a bit, how about graded geometry, graded topology, and so forth? Besides, can we define bigraded, multi-graded, $ \infty $-graded, graded of graded? 
Moreover, as the graded algebra $ universal \: enveloping \: algebra $ is universal, does ``graded'' have any connection with ``universal'' or with somewhat related ``completed''? If does, how? 
There is a rather simplified example: whether $ universal \: space \: EG $ obtained from principal $ G $-bundle is graded or not? 
If so, what actually is the graded stuff and how? $ Homotopy \: Hypothesis $ which asserts that n-groupoids are equivalent to homotopy $ n $-$ types $ for all extended natural numbers $ n\in \overline{\mathbb{N}} $ casts an illuminating model. 
Furthermore, what if the modifier graded is weakened to be stratified, sliced or ``locally structured'' in other forms? ...\\

The other one arose from a reflection on a mistake I made to take a point $ (0, g_{1}, g_{2}, ... , g_{n}) $ where $ g_{1}, g_{2}, ... , g_{n} $ are elements in topological group $ G $ in the base space of a  principal $G$-bundle  as intrinsically possessing a structure isomorphic to $ G $. 
The fact that the structure is encoded in the corresponding fiber reveals an approach to ``resolving'' singularities and ``condensed objects'' by endowing them with an extra structure attached to ``fiberlike things''. 
Setting aside techniques for resolution in singularity theory, $ Yoneda \: lemma $ along with representation theory partly annotates it, while $ sheaf \: theory $ and (topological, algebraic, equivariant, \'{e}tale, etc.) $ K$-$theory $ are known to be good supplements to fiber bundle’s. 
Therefore, it is natural to expect an appropriate method to unify these theories together. What's more, once feasible, will the unifying be possible to collide with the former idea of the ``graded''? For example, can ``loop sheaf'' connected with fundamental group be generalized to ``higher loop sheaf'' associated to higher homotopy group? 
Cohomology theory should be an easy found case that counts and deserves to be further studied under these assumptions. 
After all, homotopy theory and algebraic topology can be enlarged with a number of new elements and be excavated for profounder connotation.\\

Above thoughts, albeit superficial and naive without neither strict logic nor axiomatized semantics, propelled me to search for clues and hints in symposia, books, seminars, minicourses, and other accessible ways. 
As new ingredients from a variety of theories add to them, some questions were solved or partly solved, and new questions emerged, mixing and reacting with those partly unsolved old ones. 
The process circles over and over to make inspirations explode in my mind, driving me overwhelmed and lost in these endless puzzles. Grothendieck's theories which construct a marvelous empire with solid foundation of mathematical logic reorganized and enriched my unfettered thoughts well but don't clear my confusions on some key issues. 
Thanks to the series of mini courses on Derived Algebraic Geometry (abbreviated as $ DAG $) which was given by Bertrand To\"{e}n et cetera in October, 2019, presenting to me concrete illustration to a number of the riddles. 
Consequently, I proceeded to learn $ DAG $, starting with preliminaries like higher category theory, higher topos theory, and higher algebra, which constitute the logical basis of the renewed systematic theoretical frameworks.\\
Before diving deeper into $ DAG $, I am going to set out to explore some topics on higher orbifolds in next section. 
Since it is not finished yet despite that a plenty of work has been done, only certain parts of axiomatized semantics will be interpreted. Meanwhile, definitions, lemmas, propositions and theorems given or proved in bibliographies cited in this article would not be repeated unless necessary. 
Notwithstanding the fact that $ DAG $ theories involved in this essay are based on both Jacob Lurie's, Bertrand To\"{e}n and Gabriele Vezzosi's, for convenience, the terminology and notations involved will be borrowed from Jacob Lurie's if there are ambiguities or discrepancies between the two.\\

\section{Model Structure on Higher Orbifolds}
Orbifolds were known as singular spaces that are locally modelled on quotients of open subsets of $ R^{n} $ by finite group actions. 
The orbifold structure encodes not only the structure of the underlying quotient space, but also that of the isotropy subgroups. 
$ Atlases $ are used to describe the orbifold structure. However, it is complicated and inconvenient to define the morphisms, not to mention the composition of morphisms, hindering categorical operations on orbifolds and further steps.\\

I. Moerdijk and D. A. Pronk has shown in \cite{articleMoerdijk} that orbifolds are essentially the same as certain ``proper'' groupoids. E.  Lerman mentioned in \cite{Lerman10} that a proper \'{e}tale Lie groupoid is locally isomorphic to an finite action of groupoid. 
From then on, geometers are used to define orbifolds as proper \'{e}tale Lie groupoids \cite{moerdijk2002orbifolds}. What's more, as groupoids can also be supposed to be atlases on orbifolds, there is a way of thinking of a groupoid as “coordinates” on a corresponding stack. 
In another word, orbifolds can also be seen as stacks.\\

The concept of $ stack $ (i.e. 2-$ sheaf $) is the categorical analogue of sheaf, which is a generalization of principal bundle. 
Recall that a principal $ G $-bundle can be defined as a collection of transition functions $ g_{ij} : U_{ij} \rightarrow G $ on the double intersections $ U_{ij} $ of some open covering $ \left\lbrace U_{i} \right\rbrace _{i \in I} $ of base manifold $ M $, satisfying the cocycle condition $ g_{ij}g_{jk} = g_{ik} $ compatible with the singular cohomology group $ H^{1}(M;G) $. 
These transition functions construct morphisms of groupoids from the C\`{e}ch groupoid $ \coprod U_{ij} \rightrightarrows \coprod U_{i} $ associated to the open covering $ \left\lbrace U_{i} \right\rbrace _{i \in I} $ to the Lie groupoid $ G \rightrightarrows \ast $.  
Meanwhile, a principal $ G $-bundle P over a manifold $ M $ canonically determines a homotopy class of maps from $ M $ to the classifying space $ BG $ of the group $ G $. 
In effect, the set of isomorphism classes of $ G $-principal bundles over $ M $ is in bijection with the set of homotopy classes of maps $ M \rightarrow BG \equiv K(G,1) $, where $ K(G,1) $ denotes the Eilenberg-MacLane space. 
Thus, stacks can be viewed as a kind of classifying construction which coincides with C\`{e}ch covering.\\

Considering morphisms between  $ G $ -torsor  and $ H $ -torsor  which is generally called 2-morphisms, we can define bibundles \cite{Breen2010} ,2-groupoids and 2-stacks satisfying certain descent condition. 
In particular, $ G $-gerbes  are equivalent to $ AUT(G) $-$ principal $ 2-$ bundles $, for $ AUT(G) $ the automorphism 2-group of $ G $ \cite{43964083485245caad113a7b7d85c71a}. 
Take $ G $-bibundle as a prototype, we can construct a similar classifying stack $ BBG \equiv K(G, 2) $ associated with 2-C\`{e}ch covering satisfying 2-cocycle condition, which is compatible with $ H^{2}(M; G) $. 
As for $ H^{n}(M; G) $, where $ n \geq 2 $, $ K(G, n) $ is a higher $ n $-stack \cite{simpson1996algebraic} classifying ``higher'' pincipal $ G $-bundles \cite{articlevo}, whose $ n $-C\`{e}ch covering (essentially the same with hyper-covering) nerve is $ n $-hypergoupoid \cite{Duskin1979HigherDT}. 
Equivalently, the singular cohomology groups $ H^{n}(M; G) $ of a nice topological space $ M $ with coefficients in an abelian group $ G $ (sheaf cohomology $ H^{1}_{sheaf}(M; G) $ of $ M $ with coeﬃcients in the constant sheaf $ \mathcal{G} $ associated to $ G $ for a general space $ M $) is actually a representable functor of $ M $. 
That is, there exists an Eilenberg-MacLane space $ K(G, n) $ and a universal cohomology class $ \eta \in H ^{n} (K(G, n); G) $ such that, for any nice topological space $ X $, pullback of $ \eta $ determines a bijection  $ [X; K(G, n)] \rightarrow H^{n} (X; G) $ \cite{lurie2009higher}. 
Here $ [X; K(G, n)] $ denotes the set of homotopy classes of maps from $ X $ to $ K(G, n) $.\\

Extending $ n $ to be $ \infty $, $ \infty $-groupoid can be defined to be an $ \infty $-category \cite{lurie2009higher} (quasi-category by Joyal and weak Kan complex by Boardman and Vogt) in higher category \cite{kerodon} to be an abstract homotopical model for topological spaces. 
Therefore, stacks over $ \infty $-groupoid generalize those over groupoid, requiring a new topology corresponding to higher category, which is formally named higher topos \cite{lurie2009higher}, represented by model topos.\\

Back to orbifolds which have been previously mentioned to be proper \'{e}tale Lie groupoids. 
As Lie groupoids are essentially differentiable stacks up to Morita equivalence \cite{behrend2006differentiable}, it is reasonable to assume higher orbifolds to be “\'{e}tale differentiable higher stacks”. 
Start from higher ($ n $-)stack which can also be regarded as ($ n $-)truncated $ \infty $-stack (also denoted by ($ \infty $,n)-stack), or $ n $-geometric $ D $-stack in \cite{bertrand2008homotopical}. 
Since $ \infty $-stack (also known as $ (\infty,1) $-sheaf) is endowed with a geometry \cite{lurie2009derived}, forming a simplicial enriched category whose all internal hom-objects are $ Kan $ complexes (which are fibrant-cofibrant objects), there exists a $ n $-truncated geometry. 
We may literally transplant Jacob Lurie's frameworks in \cite{lurie2009higher} to obtain local and global model structures on $ (\infty,n) $-stacks, i.e. $ n $-orbifolds, although detailed machinery is rather complicated. 
David Carchedi gave another description in \cite{carchedi2015tale} \cite{carchedi2015homotopy} \cite{carchedi2013higher}.\\

\section{Quantum Cohomology of Higher Orbifolds}
W. Chen and Y. Ruan have asserted in \cite{chen2000orbifold} \cite{Chen_2004} that an important feature of orbifold cohomology groups is degree shifting, i.e. shifting up the degree of cohomology classes of $ X _{ ( g ) } $ by $ 2 \iota _{ ( g ) } $, defining the orbifold cohomology group of degree $ d $ to be the direct sum 
\[ 
H_{orb}^{d}(X ; \mathbb{Q})=\oplus_{(g) \in T} H^{d-2 \iota_{(g)}}\left(X_{(g)} ; \mathbb{Q}\right)
 ,\]
for any rational number $ d \in \left[ 0 , 2n \right] $, where X is a closed almost complex orbifold with $ \operatorname{dim}_{\mathbb{C}} X = n $, $ X_{(g)}=\left\{\left(p,(g)_{G_{p}}\right) \in \tilde{X} \mid(g)_{G_{p}} \in(g)\right\} $ called a twisted sector for $ (g) \neq (1) $, and $ \iota_(g) $ called degree shifting numbers. 
Moreover, Y. Ruan explored twisted orbifold cohomology and its relation to discrete torsion in \cite{ruan2000discrete}.\\

Now that the geometry and topology behavior of higher orbifolds can be detected in section 2, it should be possible to compute their quantum cohomology but will obviously be a significant amount of work. However, when it comes to derived orbifolds, things should be radically different. And what will derived twisted orbifolds look like? How would the twisted sector and degree shifting numbers change? Will them get simpler compared with the general case? Techniques are supposed to be totally different with existing literatures'.\\

I apologize to suspend it here because of time.

\bibliography{ref}
\end{document}